\documentclass[10pt,twocolumn,twoside]{IEEEtran}
\usepackage{cite}
\usepackage{epstopdf}
\usepackage{xcolor}
\usepackage[pdftex]{graphicx}
\usepackage{amsmath}
\usepackage{algorithmic}
\usepackage{array}

\usepackage{graphics} 
\usepackage{epsfig} 
\usepackage{times} 
\usepackage{amsmath} 
\usepackage{amssymb}  
\usepackage{tabularx}
\newcolumntype{b}{X}
\newcolumntype{s}{>{\hsize=.5\hsize}X}

\newcommand{\tomega}{{\tilde{\omega}}}

\newcommand{\amin}{\mathop{\mbox{argmin}}}

\newcommand{\st}{\mbox{s.t.}}

\newtheorem{prop}{Proposition}

\newtheorem{rmk}{Remark}
\newtheorem{ass}{Assumption}

\graphicspath{../figs/}

\begin{document}
\title{Computing Economic-Optimal and Stable Equilibria for Droop-Controlled Microgrids}
\author{Sungho Shin and Victor M. Zavala
  \thanks{Sungho Shin and Victor Zavala are with the Department of Chemical and Biological Engineering, University of Wisconsin-Madison, WI 53706, USA (email: {sungho.shin@wisc.edu; victor.zavala@wisc.edu}).}
}

\maketitle

\begin{abstract}
  We consider the problem of computing equilibria (steady-states) for droop-controlled, islanded, AC microgrids that are both economic-optimal and dynamically stable. This work is motivated by the observation that classical optimal power flow (OPF) formulations used for economic optimization might provide equilibria that are not reachable by low-level controllers (i.e., closed-loop unstable). This arises because OPF problems only enforce steady-state conditions and do not capture the dynamics. We explain this behavior by using a port-Hamiltonian microgrid representation. To overcome the limitations of OPF, the port-Hamiltonian representation can be exploited to derive a bilevel OPF formulation that seeks to optimize economics while enforcing stability. Unfortunately, bilevel optimization with a nonconvex inner problem is difficult to solve in general. As such, we propose an alternative approach (that we call probing OPF), which identifies an economic-optimal and stable equilibrium by probing a neighborhood of equilibria using random perturbations. {The probing OPF is advantageous in that it is formulated as a standard nonlinear program, in that it is compatible with existing OPF frameworks, and in that it is applicable to diverse microgrid models.} Experiments with the IEEE 118-bus system reveal that few probing points are required to enforce stability.
\end{abstract}


%
\IEEEpeerreviewmaketitle
\section{Introduction}
Recent advances in distributed generation (DG) technologies have enabled the development of {microgrids}, which are autonomous networks that provide flexibility to the central power grid \cite{lasseter2002microgrids,hatziargyriou2007microgrids}. Due to their desired autonomy, microgrids require robust control architectures that are capable of stabilizing them, coordinating them with neighboring networks, and maximizing their economic performance. Microgrid control is challenging due to their fast dynamics, due to disturbances covering a wide frequency spectrum (e.g., internal and external power loads and wind/solar supply), and due to the heterogeneity of the physical devices \cite{guerrero2011hierarchical,guerrero2013advanced,lopes2006defining}.

Microgrid control seeks the fundamental {\em dual} goals of maintaining maximum economic performance while also maintaining stability  \cite{dondi2002network,lopes2007integrating}. As microgrids employ various forms of controllable DG sources, the economic performance of the operation can be significantly improved by computing the optimal economic dispatch and unit commitment \cite{lasseter2002microgrids,hatziargyriou2007microgrids}. On the other hand, islanded microgrids are low-inertia power systems and thus stabilizing the operation is more challenging compared to conventional transmission/distribution systems \cite{olivares2014trends}. {\em Hierarchical} control architectures are adopted to simultaneously achieve the dual goal. The hierarchy of the microgrid control typically consists of primary, secondary, and tertiary control \cite{bidram2012hierarchical,katiraei2008microgrids,olivares2014trends}. Primary control (e.g., droop control) seeks to stabilize the system in the face of high-frequency disturbances and fast dynamic fluctuations. Secondary control (e.g., the energy management system) makes economic decisions and ensures reliable operation. Tertiary control is responsible for the coordination with the main power grid and neighboring microgrids. Since we consider islanded microgrids, we focus on primary and secondary control. In general, the reference steady-state active/reactive power generations and voltages for DG units are computed by secondary controllers, and this information is sent to primary controllers as set-points. 

Droop control is a flexible decentralized control technology used for the primary control of microgrids \cite{katiraei2005micro,gao2008control,de2004voltage,piagi2006autonomous}. Droop controllers are capable of robustly tracking the set-points given by secondary controllers. A wide range of configurations for droop control has been reported in the literature \cite{4135419,guerrero2013advanced,zhong2013robust,blaabjerg2006overview,4696040}. Despite their well-known robustness, it is difficult to analyze their closed-loop dynamic stability without making simplifying assumptions. As such, most of the existing stability analysis of droop-controlled microgrids relies on simplifying assumptions \cite{schiffer2014conditions,dorfler2013synchronization,simpson2013synchronization,dorfler2016breaking}. Schiffer et al. analyzed the stability of microgrids by assuming inductive, Kron-reduced networks \cite{schiffer2014conditions}. Such assumptions make the model less realistic but enable the derivation of powerful analytical results. Specifically, under such assumptions, droop-controlled microgrids can be modeled as {\em port-Hamiltonian} systems, for which stability can be completely characterized by a Hamiltonian function. Dorfler et al.  showed that droop-controlled microgrids are always stable under the assumptions of constant voltage magnitude (so-called {\em decoupling assumption}) and uniform resistance-to-reactance ratio \cite{dorfler2013synchronization,simpson2013synchronization,dorfler2016breaking}.

{Different strategies for secondary control have been proposed in the literature such as optimization-based approaches, expert systems (rule-based technique), and decentralized methods \cite{osti_838178}. In state-of-the-art optimization-based secondary control, off-line optimization is performed with detailed models and forecasts of renewables, demands, and market conditions \cite{palma2013microgrid} while real-time optimization is performed with current data but using simplified models \cite{huang2014adaptive}. Recently, real-time optimization-based secondary control with optimal power flow (OPF) formulations has been reported in the literature \cite{shi2017real}. OPF-based microgrid secondary control is advantageous due to its ability to directly maximize economic performance while directly enforcing network constraints \cite{4840050,6502290,dall2013distributed}. Since OPF ignores dynamics, an important question on the efficacy of this approach is whether the set-points recommended maximizing economics can be reached by the primary control layer. Recently it was shown that, for constant voltage and uniform resistance-to-reactance ratio systems, equilibria obtained with OPF are always reachable by droop control \cite{dorfler2016breaking}. However, under more general settings, the OPF formulation might provide an equilibrium that is not reachable. Such a limitation of OPF is noticed in several works (in the context of transmission networks) and motivated the development of stability-constrained OPF formulations \cite{867137,cai2008application}. Unfortunately, these approaches only consider specific types of stability criteria and do not guarantee general closed-loop stability. To the best of our knowledge, there is no general framework to ensure closed-loop stability of OPF solutions.}

In this work, we propose modifications of AC OPF formulations that compute economic-optimal and stable equilibria for general droop-controlled AC microgrids. {The framework is applicable to any type of dynamic microgrid model, network structure (meshed or radial), and economic objective (generation cost, economic revenue, or environmental factors) and can be implemented as a real-time optimization-based secondary controller.} {We first derive a dynamic model for a general droop-controlled, inverter-based microgrid that does not require any  assumptions. Such a model is general but difficult to intuitively understand its dynamic behavior. To obtain an intuitive understanding, we derive the port-Hamiltonian microgrid model by introducing additional assumptions on the network \cite{schiffer2014conditions}.} For the port-Hamiltonian model, we show that equilibrium is stable if it is a strict minimum of the Hamiltonian function. With such a property, we establish that the dynamic instability of microgrid is caused by the nonconvexity of the Hamiltonian function. Furthermore, we show that the nonconvexity arises in the region where the voltage angle differences between neighboring nodes are large. With the Hamiltonian function, we derive a bilevel OPF in which an economic objective is optimized in an outer layer and stability (expressed as a strict minimum of Hamiltonian) is enforced in an inner layer. This formulation has the key advantage that the dynamics do not need to be taken explicitly into consideration, but the nonconvexity of the Hamiltonian makes the bilevel problem computationally challenging to solve. {Moreover, the bilevel OPF cannot be applied to general microgrid models (e.g., models with non-inductive networks) because of the connection between Hamiltonian function and dynamic stability breaks if assumptions for port-Hamiltonian model do not hold.}

\begin{figure}
  \centering
  \includegraphics[width=.3\textwidth]{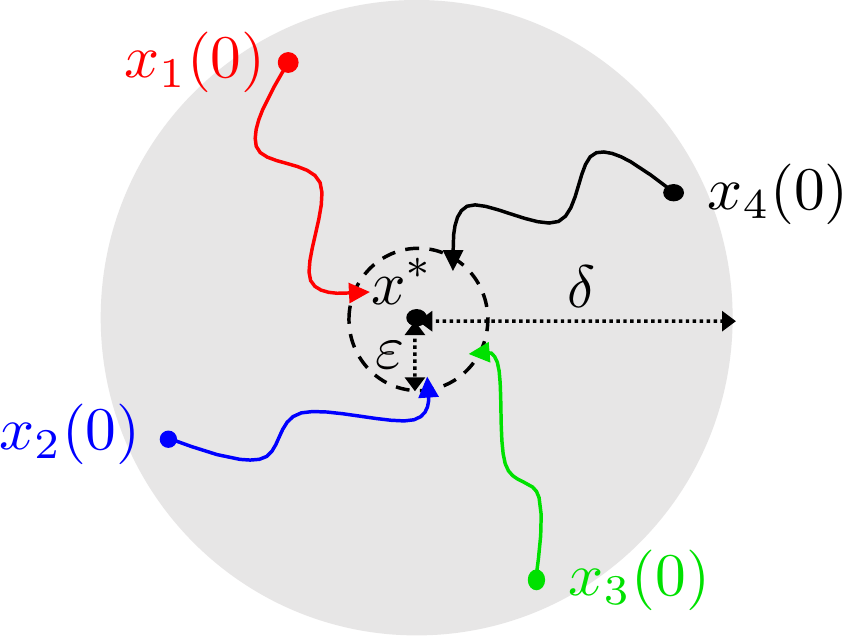}
  \caption{Schematic of probing OPF formulation, where $x^*$ represents the equilibrium, the solid lines with arrow represent the probing trajectories, the grey region represents the region of attraction with radius $\delta$, and the small dashed circle represents the allowable terminal region with radius $\varepsilon$.}\label{fig:graph-abstract}
\end{figure}

Motivated by the limitations of bilevel OPF, we propose an alternative OPF formulation that is computationally tractable. Specifically, we propose a {\em probing OPF} formulation that finds an economic-optimal and stable equilibrium by probing a neighborhood of the equilibrium using random perturbations and by enforcing the probing trajectories converge to the equilibrium (see Fig. \ref{fig:graph-abstract}). {This approach is generally applicable to different dynamic microgrid models. Also, the problem can be formulated as a single large-scale nonlinear program (NLP), whose local solution can be found in a scalable way. The probing OPF is large-scale due to the necessity to capture the system dynamics for the multiple probing trajectories, but the problem is highly sparse and structured and can be handled with existing NLP solvers \cite{Wachter2006implementation,kang2015nonlinear}. Experiments with the  IEEE 118-bus system suggest that a small number of probing points are sufficient to enforce stability.}

The paper is organized as follows. In Section \ref{sec:model}, we introduce a general mathematical model of droop-controlled, inverter-based microgrid systems. In Section \ref{sec:pHam}, we present the port-Hamiltonian model that will be used for deriving the analytical results on stability. In Section \ref{sec:convex}, the Hamiltonian function for the port-Hamiltonian model is analyzed, and we review various forms of approximate models reported in the literature. The main contributions of the paper are presented in Section \ref{sec:economic}. Here, we present the proposed bilevel and probing OPF formulations. A case study is presented in Section \ref{sec:cstudy} and conclusions are presented in Section \ref{sec:conc}.

{\em Notation.}
The set of real numbers and the set of complex numbers are denoted by $\mathbb{R}$ and $\mathbb{C}$, respectively. The real part and the imaginary part of a complex number is denoted by $\Re(\cdot)$ and $\Im(\cdot)$, respectively. We define $\mathbb{R}_{>0}:=(0,+\infty)$ and $\mathbb{R}_{\geq 0}:=[0,+\infty)$. The $i$-th component of a vector and the $(i,j)$-th component of a matrix are denoted by $(\cdot)_i$ and  $(\cdot)_{i,j}$, respectively. The transpose of a matrix or a vector is denoted by $(\cdot)^\top$.
  We consider all vectors as column vectors and use the syntax $(v_1,v_2,\cdots ,v_n) = [v_1^\top v_2^\top \cdots v_n^\top ]^\top$. To represent vectors, we use the notation $\{x_i\}_{i\in\mathcal{I}}:=(x_{i_1},x_{i_2},\cdots,x_{i_n})$, where $\mathcal{I}:=\{i_1<i_2<\cdots<i_n\}$. Furthermore, we write $x_\mathcal{I}:=\{ x_i \}_{i\in\mathcal{I}}$. The $n\times n$ identity matrix is denoted by $I_n$, and the subscript is suppressed if the dimension can be implicitly deduced. The power networks discussed in this work are undirected weighted graphs $\mathcal{G}(\mathcal{V},\mathcal{E},w)$ where $\mathcal{V}$ is the set of nodes (buses) and $\mathcal{E}$ is the set of edges (lines). The undirected edge between node $i$ and $j$ is denoted by $\{i,j\}$. The weight of edge $\{i,j\}$ is denoted by $(w)_{i,j}=(w)_{j,i}$. When variables $x_i$ are defined for every node $i\in\mathcal{V}$, we use the compact notation $x:=x_\mathcal{V}$. Open and closed neighborhoods of node $i\in\mathcal{V}$ are denoted by $\mathcal{N} (i)$ and $\mathcal{N}[i]$, respectively, and defined as $\mathcal{N}(i):=\{j\in\mathcal{V}\,\mid \{i,j\}\in\mathcal{E}\}$ and $\mathcal{N}[i]:=\mathcal{N}(i)\cup\{i\}$, respectively. A flow $x(\cdot)$ generated by a dynamical system with an initial value $x^0$ is denoted by $x(\cdot;x^0)$. We denote cardinalities of sets and absolute values of real numbers by $\vert\cdot\vert$. Euclidean vector norms and induced Euclidean norms of matrices are denoted as $\Vert\cdot\Vert$.

\section{Droop-Controlled Microgrid Dynamics}\label{sec:model}
In this section, we derive the dynamics of a droop-controlled, inverter-based microgrid model.
The detailed setting is as follows. The network is described as a connected, undirected, complex-weighted graph $\mathcal{G}(\mathcal{V},\mathcal{E},w)$ with the node set $\mathcal{V}:=\{1,2,\cdots,n\}$. The weight $(w)_{i,j}\in\mathbb{C}$ of edge $\{i,j\}\in\mathcal{E}$ represents the admittance of the edge $\{i,j\}$. The admittance matrix $Y\in\mathbb{C}^{n\times n}$ is defined as the weighted Laplacian of the graph. The set of generator nodes are denoted by $\mathcal{V}_G:=\{1,\cdots,m\}$ and the set of load nodes are denoted by $\mathcal{V}_L:=\{m+1,\cdots,n\}$. We assume that each generator node is composed of inverters interfaced with DG sources (e.g. photovoltaic cells, microturbines, batteries). The microgrid dynamics can be described by the trajectories of the voltage magnitudes $V_i:\mathbb{R}_{\geq 0} \rightarrow \mathbb{R}_{>0}$ and angles and $\theta_i:\mathbb{R}_{\geq 0} \rightarrow \mathbb{R}$ for each node $i\in\mathcal{V}$. Lastly, we consider the last generator node $m$ as a reference node (i.e., $\theta_m\equiv 0$). The dynamic model comprises power flow equations, droop-control laws, filter dynamics, and load models. 

{ 
\subsection{Power Flow Equations}\label{sec:pflow}
Active and reactive power injections $P_{i}:\mathbb{R}^{\vert \mathcal{N}[i]\vert }\times\mathbb{R}_{>0}^{\vert \mathcal{N}[i]\vert } \rightarrow \mathbb{R}$ and $Q_{i}:\mathbb{R}^{\vert \mathcal{N}[i]\vert }\times\mathbb{R}_{>0}^{\vert \mathcal{N}[i]\vert } \rightarrow \mathbb{R}$ of node $i\in\mathcal{V}$ are functions of the voltage angles $\theta_{\mathcal{N}[i]}$ and magnitudes $V_{\mathcal{N}[i]}$ of their neighboring (connected) nodes and can be expressed as
\begin{subequations}\label{eqn:pflow-gen}
  \begin{align}
    \label{eqn:apflow-gen}
    P_{i} &= \sum_{j\in\mathcal{N}[i]}V_i V_j\left( G_{ij} \cos\theta_{ij} + B_{ij}  \sin \theta_{ij}\right), \quad i\in\mathcal{V}\\
    \label{eqn:rpflow-gen}
    Q_{i} &= \sum_{j\in\mathcal{N}[i]} V_i V_j\left(G_{ij} \sin\theta_{ij} - B_{ij} \cos \theta_{ij}\right), \quad i\in\mathcal{V},
  \end{align}
\end{subequations}
where $G_{ij}=G_{ji}:=\Re(Y_{ij})$, $B_{ij}=B_{ji}:=\Im(Y_{ij})$, and $\theta_{ij}:=\theta_i-\theta_j$. Note that the dependence of the active and reactive power on $\theta_{\mathcal{N}[i]}$ and $V_{\mathcal{N}[i]}$ is suppressed from the notation. }

\subsection{Droop-Control Laws and Filter Dynamics}
We model the inverters as AC voltage sources for which voltage magnitudes and frequencies can be manipulated.
\begin{subequations}\label{eqn:droop-law}
\begin{align}\label{eqn:inv}
  \omega_i &= u_i^\omega,\quad V_i = u_i^V,\quad i\in\mathcal{V}_G
\end{align}
Here, $\omega_i:=\omega_m+d\theta_i/dt$ is the frequency (note that $\theta_i$ is relative voltage angle in the frame of node $n$) and $u_i^\omega:\mathbb{R}_{\geq 0} \rightarrow \mathbb{R}$ and $u_i^V:\mathbb{R}_{\geq 0} \rightarrow \mathbb{R}$ are the control variables for the inverters of node $i$. We use the conventional active power-frequency ($P-\omega$) and reactive power-voltage magnitude ($Q-V$) droop \cite{chandorkar1993control}. In particular, the droop-controller of node $i$ manipulates the voltage frequency and magnitude of node $i$ using the measured active and reactive power injections of node $i$ using the following control laws:
\begin{align}\label{eqn:droop-p}
  u_i^\theta &= \omega^d - k_{P_i} (P_i^m - P^d_i),\quad i\in\mathcal{V}_G\\
  u_i^V &= V^d_i - k_{Q_i} (Q_i^m - Q^d_i),\quad i\in\mathcal{V}_G.\label{eqn:droop-q}
\end{align}
Here, $P_i^m:\mathbb{R}_{\geq 0} \rightarrow \mathbb{R}$ and $Q_i^m:\mathbb{R}_{\geq 0} \rightarrow \mathbb{R}$ are the measured active and reactive power injections of node $i$, $\omega^d, V_i^d, P_i^d, Q_i^d\in\mathbb{R}$ are the {\em set-points} of the voltage frequency and magnitude, and the active and reactive power injections of node $i\in\mathcal{V}$, and $k_{P_i},\ k_{Q_i}\in\mathbb{R}_{>0}$ are the droop control gains. The control gains are typically determined from the steady-state characteristics of the system \cite{chandorkar1993control,li2004design,katiraei2006power}.

It is assumed the active and reactive powers are measured with first-order filters as
\begin{align}\label{eqn:fltr}
  \tau_{P_i} \dot{P}_i^m(t) = P_i- P_i^m,\; \tau_{Q_i} \dot{Q}_i^m(t) = Q_i- Q_i^m,
\end{align}
\end{subequations}
for $i\in\mathcal{V}_G$, where $\tau_{P_i}$, $\tau_{Q_i}\in\mathbb{R}_{>0}$ are the filter constants. {Note that model \eqref{eqn:droop-law} can be considered as a reduced-order inverter-based microgrid model that can be formally obtained from the singular perturbation analysis of the fundamental physics-based model \cite{schiffer2016survey,ajala2018hierarchy,luo2014spatiotemporal}. The first-order dynamic model that has been used in a number of other works in the literature can be obtained in the special case where $\tau_{P_i}=\tau_{Q_i}=0$ for $i\in \mathcal{V}_G$.}

{\begin{rmk}[Droop controller configuration]\label{rmk:other-models}
  We use a $P-\omega$ and $Q-V$ droop control scheme \eqref{eqn:droop-law}, but many other configurations are possible (see \cite{4135419,guerrero2013advanced,zhong2013robust}). Furthermore, there can be many other droop control settings other than the PI control scheme used in \eqref{eqn:droop-law}, ranging from simple P-control to PI control with feedforward compensation (see \cite{olivares2014trends,blaabjerg2006overview,4696040}). We highlight that the proposed probing OPF formulation, presented in Section \ref{sec:economic}, does not require a specific algebraic form of the model; as long as the model can be expressed as a system of differential and algebraic equations (DAEs). Here, we use a specific example of droop control setting \eqref{eqn:droop-law} to describe the typical structure of dynamic microgrid models and for the sake of consistency with the port-Hamiltonian model discussed in existing literature studies.
\end{rmk}}
{\subsection{Load Models}
We use the standard constant active and reactive power load model for each load bus (i.e., the loads are not dynamic)
\begin{align}\label{eqn:load}
  P^d_i = P_i,\quad Q^d_i=Q_i,\quad i\in\mathcal{V}_L,
\end{align}
where $P_i^d, Q_i^d\in\mathbb{R}$ are the active and reactive power loads. Note that we can allow generator nodes to have non-zero loads. In such a case, we consider $P^d_i$ and $Q^d_i$ for $i\in\mathcal{V}_G$ as the set-points for the net active and reactive power injections (set-point for generation$-$load).}

\subsection{Full Closed-Loop Model}
The general droop-controlled microgrid model can be derived by simplifying \eqref{eqn:pflow-gen}-\eqref{eqn:load}. By differentiating \eqref{eqn:inv}-\eqref{eqn:droop-q} and eliminating $u^\omega_i,\ u^V_i, P_i^m$, and $ Q_i^m$, we obtain
\begin{subequations}\label{eqn:ss-gen}
  \begin{align}
    \dot{\theta_i}(t) &= \omega_i-\omega_n
    ,\quad i\in\mathcal{V}_G\setminus\{m\}\\
    \tau_{P_i} \dot{\omega}_i(t)&= -\omega_i + \omega^d - k_{P_i} (P_i - P_i^d)\label{eqn:freq}
    ,\quad i\in\mathcal{V}_G\\
    \tau_{Q_i} \dot{V}_i(t)&= -V_i + V_i^d - k_{Q_i} (Q_i - Q_i^d)
    ,\quad i\in\mathcal{V}_G\\
     P^d_i &= P_i ,\quad Q^d_i=Q_i \label{eqn:pflow-load}
    ,\quad i\in\mathcal{V}_L.
  \end{align}
\end{subequations}
We refer to \eqref{eqn:ss-gen} as a {\em general model}, as opposed to the port-Hamiltonian model discussed later. Observe that the voltage angles are defined for all $\mathcal{V}\setminus \{m\}$ (node $m$ is considered as a reference node), the voltage magnitudes are defined for all $i\in\mathcal{V}$, and the voltage frequencies are defined for all $\mathcal{V}_G$. Thus, we denote the full state vector as $x:=(\theta,\omega,V)\in\mathbb{R}^{n-1}\times \mathbb{R}^{m}\times \mathbb{R}_{>0}^n$. We consider $\omega^d$, $V^d_i$,$P^d_i$, and $Q^d$ for $i\in\mathcal{V}_G$ as inputs provided by the secondary control layer. Furthermore, the loads $P_i^d$ and $Q_i^d$ for $i\in\mathcal{V}_L$ can be regarded as uncontrollable inputs. By inspecting \eqref{eqn:ss-gen}, we observe that it suffices to consider the following as inputs for $i\in\mathcal{V}$.
\begin{align}\label{eqn:u}
  P^u_i &: =
  \begin{cases}
    P_i^d+\frac{\omega^d}{k_{P_i}}& i\in\mathcal{V}_G\\
    P_i^d& i\in\mathcal{V}_L
  \end{cases}\;
  Q^u_i := 
  \begin{cases}
    Q_i^d+\frac{V^d_i}{k_{Q_i}}&i\in\mathcal{V}_G\\
    Q_i^d&i\in\mathcal{V}_L
  \end{cases}
\end{align}
The full input vector is denoted as $ u:=(P^u, Q^u)\in\mathbb{R}^{2n}$. 

{The droop-controlled microgrid model \eqref{eqn:ss-gen} is a set of semi-explicit DAEs. The model can be written in the following form
  \begin{align}\label{eqn:gen-dae}
    (\dot{\theta}_{\mathcal{V}_G\setminus\{m\}},\dot{\omega}, \dot{V}_{\mathcal{V}_G}) = g(x,u),\quad 0=h(x,u)
  \end{align}
  with $g:(\mathbb{R}^{n-1}\times\mathbb{R}^{m}\times\mathbb{R}_{>0}^{n})\times \mathbb{R}^{2n}\rightarrow \mathbb{R}^{3m-1}$ and $h:(\mathbb{R}^{n-1}\times\mathbb{R}^{m}\times\mathbb{R}_{>0}^{n})\times \mathbb{R}^{2n}\rightarrow \mathbb{R}^{2(n-m)}$. We assume that the algebraic equations in \eqref{eqn:gen-dae} are always solvable for given $(\theta_{\mathcal{V}_G\setminus\{m\}},V_{\mathcal{V}_G})$ and $u$. That is, the semi-explicit DAE \eqref{eqn:ss-gen} is of index one. We can thus write \eqref{eqn:ss-gen} as an ordinary differntial equation (ODE):
\begin{align}\label{eqn:gen}
  \dot{x}&=f(x,u).
\end{align}
Here we assume that $f:(\mathbb{R}^{n-1}\times\mathbb{R}^{m}\times\mathbb{R}_{>0}^{n})\times \mathbb{R}^{2n}\rightarrow \mathbb{R}^{2n+m-1}$ is twice continuously differentiable. For the sake of simplicity, we interchangably use either \eqref{eqn:gen-dae} or \eqref{eqn:gen} to represent the dynamic moodel in a compact form. In Section \ref{sec:economic}, we show that the steady-state condtitions ($f(x,u)=0$ or equivalently, $g(x,u)=0$ and $h(x,u)=0$) are satisfied at the solution of the OPF problem. 

\begin{rmk}[Stability condition of \eqref{eqn:ss-gen}]\label{rmk:Hurwitz}
  Consider the steady-state $(x^s,u^s)$ of the dynamic microgrid model \eqref{eqn:gen} (i.e., $f(x^*,u^*)=0$). The dynamic stability condition of \eqref{eqn:gen} states that the Jacobian $\frac{\partial f}{\partial x}(x^*,u^*)$ of $f(\cdot,u^*)$ evaluated at $x^*$ is Hurwitz \cite{khalil1996nonlinear}. That is, all the eigenvalues $\lambda$ of the Jacobian matrix satisfy $\Re(\lambda)<0$. There is no guarantee that the Jacobian is Hurwitz. Moreover, it is difficult to establish analytical sufficient conditions for that without making simplifying assumptions. In Section \ref{sec:pHam}, we discuss a simplified model and establish sufficient conditions for stability. In Section \ref{sec:convex}, we demonstrate that the stability can be indeed violated. Thus, merely finding a steady-state does not guarantee stability.
\end{rmk}
}

\section{Port-Hamiltonian Microgrid Model}\label{sec:pHam}
{In this section, we derive the port-Hamiltonian microgrid model suggested by Schiffer et al. \cite{schiffer2014conditions}. With the simplified model, we can show that the microgrid model might indeed deliver unstable equilibria. We emphasize that we discuss the simplified models only to highlight the potential instability issue of microgrids, and those simplifying assumptions are not necessary for the probing OPF.

\subsection{Model Derivation}
\begin{ass}[Port-Hamiltonian model asssumptions]\label{ass:pHam}
  Consider a droop-controlled, inverter-based microgrid model \eqref{eqn:ss-gen}. We make the following assumptions.
  \begin{enumerate}
  \item [(a)] The network is purely inductive and thus $G_{ij}=0 ,\quad \forall i,j\in\mathcal{V}$. 
  \item [(b)] The network only consists of generator nodes and thus $\mathcal{V}_G=\mathcal{V},\quad \mathcal{V}_L=\emptyset$. 
  \end{enumerate}
\end{ass}
Assumption \ref{ass:pHam} considers an idealized microgrid model, and this can be unrealistic in practical settings. The following remarks discuss circumstances under which Assumption \ref{ass:pHam} can be justified.
\begin{rmk}
  Conventionally, the lossless assumption (Assumption \ref{ass:pHam}(a)) was prevalent in power flow studies that target high-voltage transmission networks. However, in low and medium-voltage power networks like microgrids, the network can be more resistive than inductive. The lossless assumption holds only in special cases where the network is dominantly inductive due to the effect of primarily inductive inverter outputs and/or due to the presence of transformers \cite{schiffer2014conditions}.
\end{rmk}
\begin{rmk}
  Assumption \ref{ass:pHam}(b) can be justified by assuming that Kron reduction \cite{kundur1994power} is performed to eliminate the algebraic equations associated with loads. In this case, the network $\mathcal{G}(\mathcal{V},\mathcal{E},w)$ is a Kron-reduced network, only composed of generator nodes (each generator nodes are allowed to have loads). Alternatively, one could consider the loads as energy storage-interfaced inverters or dynamic loads (this is called a {\em structure-preserving model} \cite{622983,5208252}). 
\end{rmk}}

Under Assumption \ref{ass:pHam}, the model reduces to
\begin{subequations}\label{eqn:ss0}
  \begin{align}
    \dot{\theta_i}(t) &= \omega_i - \omega_n
    ,&&\; i\in\mathcal{V}\setminus\{n\}\\
    \tau_{P_i} \dot{\omega}_i(t)&= -\omega_i + \omega^d - k_{P_i} (P_i - P_i^d)\label{eqn:freq0}
    ,&&\; i\in\mathcal{V}\\
    \tau_{Q_i} \dot{V}_i(t)&= -V_i + V_i^d - k_{Q_i} (Q_i - Q_i^d)
    ,&&\; i\in\mathcal{V},
  \end{align}
\end{subequations}
and the power flow equations reduce to 
\begin{align*}
  P_{i} &= \sum_{j\in\mathcal{N}(i)}V_i V_jB_{ij}  \sin \theta_{ij},\;
  Q_{i} = \sum_{j\in\mathcal{N}[i]}  -V_i V_j B_{ij} \cos \theta_{ij},
\end{align*}
for $i\in\mathcal{V}$. Note that now $n$ is the reference node.

For convenience, we redefine voltage frequencies $\omega_i$ as relative quantities. Suppose that all the nodes are synchronized to a common frequency (i.e., $\omega_1=\cdots=\omega_n$ and $\dot{\omega}_1(t)=\cdots=\dot{\omega}_n(t)=0$). Such a frequency can be calculated by dividing \eqref{eqn:freq0} by $k_{P_i}$ and summing over all the nodes $i\in\mathcal{V}$. By inspecting the power flow equations, we can observe that $\sum_{i\in\mathcal{V}}P_i=0$ (a property of the lossless network). Thus, the synchronized frequency $\omega^s$ can be expressed as:
\begin{align*}
  \omega^s := \omega^d - \frac{\sum_{i\in\mathcal{V}} P_i^d }{\sum_{i\in\mathcal{V}} 1/k_{P_i}}.
\end{align*}
We define the new variables $\tomega_i := \omega_i - \omega^s,\, i\in\mathcal{V}$ and $\tomega^d := \omega^d - \omega^s$. By using some abuse in terminology, we refer to $\tomega_i$ as the voltage frequency of node $i$. Using this new notation, system  \eqref{eqn:ss0} reads as follows.
\begin{subequations}\label{eqn:ss}
  \begin{align}
    \dot{\theta_i} &= \tomega_i - \tomega_n
    ,&& i\in\mathcal{V}\setminus \{n\}\label{eqn:sst}\\
    \tau_{P_i} \dot{\tomega}_i&= -\tomega_i + \tomega^d - k_{P_i} (P_i - P_i^d)
    ,&& i\in\mathcal{V}\label{eqn:sso}\\
    \tau_{Q_i} \dot{V}_i&= -V_i + V_i^d - k_{Q_i} (Q_i - Q_i^d)
    ,&& i\in\mathcal{V}\label{eqn:ssV}
  \end{align}
\end{subequations}
We refer to \eqref{eqn:ss} as port-Hamiltonian model.
We denote the full state vector as $x:=(\theta,\tomega,V)\in\mathbb{R}^{2n-1}\times \mathbb{R}_{>0}^n$. For port-Hamiltonian model, we redefine $P^u_i:=P^d_i+\frac{\tilde{\omega}^d}{k_{P_i}}$ and $Q^u_i$ is defined the same as in \eqref{eqn:u}. Lastly, we assume that the input is always chosen in a way that satisfies $\tomega^d=0$ (i.e., $\sum_{i=1}^n P^u_i=0$). This constraint on $u$ can be incorporated in the secondary control layer. Observe that all the algebraic equtations associated with the loads are now replaced by the differential equations and the model reduces to a simple ODE form \eqref{eqn:gen} without any modification. 

\subsection{Hamiltonian Function of Closed-Loop System}\label{sec:convexiltonian}
System \eqref{eqn:ss} is a port-Hamiltonian system; the dynamic mapping $f(\cdot,u)$ can be expressed  in terms of the gradient of the Hamiltonian function $H:(\mathbb{R}^{2n-1}\times\mathbb{R}_{>0}^{n})\times \mathbb{R}^{2n}\rightarrow\mathbb{R}$ defined as follows.
\begin{align}\label{eqn:ham}
  H(x,u):=&\sum^n_{i\in\mathcal{V}} \left( \frac{\tau_{P_i}{\tilde{\omega}}_i^2}{2k_{P_i}} 
  + \frac{V_i}{k_{Q_i}} -Q^u_i \ln V_i
  -\frac{B_{ii} V_i^2}{2} - P^u_i\theta_i\right)\nonumber \\
  &- \sum_{\{i,j\}\in\mathcal{E}} B_{ij} V_i V_j \cos \theta_{ij}
\end{align}
By inspecting \eqref{eqn:ham}, we see that $H(\cdot,u)$ is a twice continuously differentiable mapping with respect to $x$ on $\mathbb{R}^{2n-1}\times\mathbb{R}^n_{>0}$ and for any $u\in\mathbb{R}^{2n}$. The partial derivatives are
\begin{subequations}\label{eqn:pdev}
  \begin{align}
    \frac{\partial H}{\partial \theta_i}&=-P^u_i + P_i\\
    \frac{\partial H}{\partial \tomega_i}&=\frac{\tau_{P_i}\tomega_i}{k_{P_i}}\\
    \frac{\partial H}{\partial V_i}&=\frac{1}{k_{Q_i}}-\frac{Q^u_i-Q_i}{V_i} .
  \end{align}
\end{subequations}
The mappings $f(\cdot,u)$ and $H(\cdot,u)$ are related as
\begin{align*}
  f(x,u) = (J-R(x)) \nabla_x H(x,u),
\end{align*}
with
\begin{align*}
  R(x) &=
  \begin{bmatrix}
    &&\\
    & R_{\tomega\tomega} & \\
    &  & R_{V V}(x)\\
  \end{bmatrix}\;
  J =
  \begin{bmatrix}
    & J_{\theta\tomega} & \\
    -J_{\theta\tomega}^\top &  & \\
    &  & \\
  \end{bmatrix},
\end{align*}
and $R_{\tomega\tomega},R_{VV}\in\mathbb{R}^{n\times n}$, $J_{\theta\tomega}\in\mathbb{R}^{(n-1)\times n}$ are defined by:
\begin{subequations}\label{eqn:RJs}
  \begin{align}
    (R_{\tomega\tomega})_{ij}&:=
    \begin{cases}
      \frac{k_{P_i}}{\tau_{P_i}^2} & i=j\\
      0 & i\neq j
    \end{cases}
    \quad
    (R_{VV}(x))_{ij}:=
    \begin{cases}
      \frac{k_{Q_i}}{\tau_{Q_i}}V_i & i=j\\
      0 & i\neq j
    \end{cases}\\
    (J_{\theta\tomega})_{ij}&:=
    \begin{cases}
      \frac{k_{P_i}}{\tau_{P_i}} & i=j\\
      -\frac{k_{P_i}}{\tau_{P_i}} & j=n\\
      0 & \text{otherwise}.
    \end{cases}
  \end{align}
\end{subequations}
Here, we used $\sum_{i\in\mathcal{V}} P^u_i=\sum_{i\in\mathcal{V}} P_i=0$. We note that, for all $x\in\mathbb{R}^{2n-1}\times\mathbb{R}_{>0}^{n}$, $R(x)\in\mathbb{R}^{(3n-1)\times(3n-1)}$ is positive semi-definite (PSD), $J\in\mathbb{R}^{(3n-1)\times(3n-1)}$ is skew-symmetric, and $J-R(x)$ is nonsingular.  The nonsingularity of $J-R(x)$ can be established by obtaining the Schur complement \cite[pg 22]{horn1990matrix} on $R_{\tomega\tomega}$ and noting that $R_{VV}(x)$ and $J_{\theta\tomega} R_{\tomega\tomega}^{-1} J_{\theta\tomega}^\top $ are nonsingular. 

A key property is that, for any input $u\in\mathbb{R}^{2n}$, the Hamiltonian decreases monotonically along a trajectory $x(\cdot;x_0)$ with any starting point $x^0\in\mathbb{R}^{2n-1}\times\mathbb{R}_{>0}^{n}$. One can show
\begin{align*}
  \dot{H}(x,u)
  &= (\nabla H)^\top \dot{x} = (\nabla H)^\top f(x,u)\\
  &=-\nabla_x H(x,u)^\top R(x) \nabla_x H(x,u),
\end{align*}
and thus,
\begin{align}\label{eqn:mono}
  \dot{H}(x,u)   \leq 0.
\end{align}
Inequality \eqref{eqn:mono} implies that the Hamiltonian can be used as a Lyapunov function of the closed-loop system.
\subsection{Set-Point Reachability}
The reachability of a given set-point (equilibrium) can be analyzed by examining the closed-loop asymptotic stability \cite{khalil1996nonlinear}. Given $u$, an equilibrium $x^*\in \mathbb{X}$ of a system $\dot{x}=f(x,u)$ is asymptotically stable on $\mathbb{X}$ if $\lim_{t\rightarrow\infty} x(t;x_0) = x^*$ for any $x_0\in\mathbb{X}$.
The following proposition generalizes the results in Schiffer et al. \cite{schiffer2014conditions} and provides a sufficient condition for asymptotic stability.

\begin{prop}\label{prop:stab}
  Suppose that system $\dot{x}=f(x,u)$ on a domain $\mathbb{D}\subseteq\mathbb{R}^{n_x}$ with non-empty interior has a twice continuously differentiable Hamiltonian $H(\cdot,u)$ (on $\mathbb{D}$ with respect to $x$) satisfying:
  \begin{subequations}\label{eqn:resilience}
    \begin{align}
      \dot{H}(x,u) &\leq 0\label{eqn:resilience1}\\
      \dot{H}(x,u) &= 0 \implies \nabla_{x} H(x,u) = 0\label{eqn:resilience2}  
    \end{align}
  \end{subequations}
  for any $x\in \mathbb{D}$. Then, any local strict minimum $x^*$ of $H(\cdot,u)$ in $\mathbb{D}$ with $\nabla_{xx}H(x^*,u)>0$ is locally asymptotically stable.
\end{prop}
\begin{IEEEproof}
  The proof is given in Appendix \ref{sec:pf-stab}.  
\end{IEEEproof}

The droop-controlled microgrid \eqref{eqn:ss} satisfies assumptions of Proposition \ref{prop:stab} with $\mathbb{D}=\mathbb{R}^{2n-1}\times\mathbb{R}^n_{>0}$; the Hamiltonian  $H(\cdot,u)$ is twice continuously differentiable, the monotonic decrease of the Hamiltonian (shown in \eqref{eqn:mono}) satisfies \eqref{eqn:resilience1}, and the nonsingularity of $J-R(x)$ yields \eqref{eqn:resilience2}. In most applications, strict minimizers satisfy the second order sufficient solutions. We thus have that strict minima of the Hamiltonian are equilibra that are reachable by droop control. 

\section{Convexity Analysis of Hamiltonian}\label{sec:convex}
We now investigate conditions that can guarantee the convexity of the Hamiltonian \eqref{eqn:ham}.  This is important because, if the Hamiltonian is strictly convex, any equilibrium point of \eqref{eqn:ss} is asymptotically stable. On the other hand, if the Hamiltonian is nonconvex, the steady-state condition cannot guarantee stability. We also extend the convexity analysis results to other existing approximate microgrid models \cite{dorfler2013synchronization,simpson2013synchronization,dorfler2016breaking}.

\subsection{The Hessian of the Hamiltonian}
Convexity can be assessed by analyzing the Hessian of the Hamiltonian. By taking partial derivatives of \eqref{eqn:pdev}, we obtain
\begin{align*}
  \nabla_{xx}H(x,u) =
  \begin{bmatrix}
    L(x) & & W(x)\\
    & A & \\
    W(x)^\top & & D(u)+T(x)
  \end{bmatrix}
\end{align*}
where $L(x)\in\mathbb{R}^{(n-1)\times(n-1)}$, $W(x)\in\mathbb{R}^{(n-1)\times n} $, $T(x)\in\mathbb{R}^{n\times n} $ are defined as follows:
\begin{align*}
  (L(x))_{i,j} &:= 
  \begin{cases}
    \sum_{k\in\mathcal{N}(i)} B_{ik}V_i V_k \cos\theta_{ik} &\text{if}\quad i=j\\
    -B_{ij}V_i V_j \cos\theta_{ij} &\text{if}\quad i\neq j
  \end{cases}\\
  (W(x))_{i,j} &:= 
  \begin{cases}
    \sum_{k\in\mathcal{N}(i)} B_{ik} V_k \sin\theta_{ik} &\text{if}\quad i=j\\
    B_{ij}V_i \sin\theta_{ij} &\text{if}\quad i\neq j
  \end{cases}\\
  (T(x))_{i,j} &:= 
  \begin{cases}
    -B_{ii} &\text{if}\quad i=j\\
    -B_{ij} \cos\theta_{ij} &\text{if}\quad i\neq j,
  \end{cases}
\end{align*}
and $A,D(u)\in\mathbb{R}^{n\times n}$ are defined as follows
\begin{subequations}\label{eqn:AD}
  \begin{align}\label{eqn:A}
    (A)_{i,j} &:= 
    \begin{cases}
      \frac{\tau_{P_i}}{k_{P_i}}  &\text{if}\quad i=j\\
      0 &\text{if}\quad i\neq j
    \end{cases}
    \\
    (D(u))_{i,j} &:= 
    \begin{cases}
      \frac{Q^u_i}{V_i^2}  &\text{if}\quad i=j\\
      0 &\text{if}\quad i\neq j\label{eqn:D}.
    \end{cases}
  \end{align}
\end{subequations}
We now show that, under mild assumptions, there exists a region where $\nabla_{xx}H(x,u)$ is positive definite (PD). Observe that the Hessian is PD if and only if the Schur complement
\begin{align}\label{eqn:S}
  S(x,u) := L(x) - W(x)^\top (D(u)+T(x))^{-1} W(x)  
\end{align}
is PD. We can see that $L(x)$ is PD if $\theta\in\Theta_G(\pi/2)$ but $W(x)^\top (D(u)+T(x))^{-1} W(x)$ can make $S(x,u)$ indefinite, where
\begin{align}\label{eqn:thg}
  \Theta_G(\gamma):=\{\theta\in\mathbb{R}^{n-1}\mid |\theta_{ij}|<\gamma \ \text{for} \ \{i,j\}\in\mathcal{E}\}.
\end{align}
The following establishes the region where $\nabla_{xx}H(x,u)>0$.
\begin{prop}\label{prop:convex}
Consider the port-Hamiltonian microgrid model \eqref{eqn:ss} with Assumption \ref{ass:pHam}. Suppose that $Q^d_i+V^d_i/k_{Q_i}>0$ and $V_{\min}<V_i<V_{\max}$ for $i\in\mathcal{V}$, then there exists $\epsilon\in\mathbb{R}_{>0}$ such that if $\theta\in\Theta_G(\epsilon)$, then $ \nabla_{xx} H(x,u) >0$.
\end{prop}
\begin{IEEEproof}
The proof is given in Appendix \ref{sec:pf-convex}.  
\end{IEEEproof}

Proposition \ref{prop:convex} indicates that, if the voltage angle differences between the neighboring nodes are sufficiently small, the Hessian is PD, and thus the Hamiltonian is strictly convex.

\begin{rmk}[Cause of Instability]
  Proposition \ref{prop:convex} reveals conditions under which an equilibrium of droop-controlled microgrid is guaranteed to be stable. Conversely, if such conditions are not satisfied, there is no guarantee for stability. The analysis also reveals that, when the voltage angle separations are not tightly bounded, the state can potentially enter the region where Hamiltonian is nonconvex. This leads to an observation that OPF is likely to find a stable equilibrium if tight bound constraints are enforced to voltage angle separations. However, such a conservative constraint setting can cause the loss of economic performance and the infeasibility of the problem. Thus, stability should be {\em directly enforced} by different means.
\end{rmk}

\subsection{Decoupling Approximation}
Now we review a different approach for establishing the stability condition of droop-controlled microgrids. Dorfler et al. \cite{dorfler2013synchronization,simpson2013synchronization,dorfler2016breaking} apply a decoupling approximation where voltage magnitudes are assumed to be constant (i.e., dynamics caused by reactive power are disregarded). Under such an assumption, the model reduces to:
\begin{subequations}\label{eqn:ss-dorfler}
  \begin{align}
    \dot{\theta_i} &= \tomega_i - \tomega_n
    ,\quad i\in\mathcal{V}\setminus \{n\}\\
    \tau_{P_i} \dot{\tomega}_i&= -\tomega_i + \tomega^d - k_{P_i} (P_i - P_i^d)
    ,\quad i\in\mathcal{V}\\
    P_{i} &= \sum_{j\in\mathcal{N}(i)} B_{ij} \tilde{V}_i \tilde{V}_j\sin \theta_{ij}, \quad i\in\mathcal{V}
  \end{align}
\end{subequations}
where $\tilde{V}_i\in\mathbb{R}_{>0}$ are fixed. We let $\tilde{x}:=(\theta,\tomega)$, $\tilde{u}:=P^u$, and
\begin{align*}
  \tilde{H}(\tilde{x},\tilde{u}):=&\sum^n_{i=1} \left( \frac{\tau_{P_i}{\tilde{\omega}}_i^2}{2k_{P_i}} - P^u_i\theta_i\right)
  - \sum_{\{i,j\}\in\mathcal{E}} B_{ij}\tilde{V}_i \tilde{V}_j \cos \theta_{ij}.
\end{align*}
The dynamics become $\dot{\tilde{x}} = (\tilde{R}-\tilde{J})\nabla_{\tilde{x}} \tilde{H}(\tilde{x},\tilde{u})$,  where:
\begin{align*}
  \tilde{R} &:=
  \begin{bmatrix}
    &  \\
    & R_{\tomega\tomega} 
  \end{bmatrix}
  \quad
  \tilde{J} :=
  \begin{bmatrix}
    & J_{\theta\tomega} \\
    -J^\top_{\theta\tomega} &  
  \end{bmatrix}
\end{align*}
and $R_{\tomega\tomega}$ and $J_{\theta\tomega}$ are defined in \eqref{eqn:RJs}. We have that $\tilde{R}$ is PSD, $\tilde{J}$ is skew-symmetric, and $\tilde{J}-\tilde{R}$ is nonsingular. Thus, \eqref{eqn:resilience} is satisfied and Proposition \ref{prop:stab} can be applied (i.e., the minimizer of $\tilde{H}(\cdot,\tilde{u})$ is a stable equilibrium of \eqref{eqn:ss-dorfler}). A similar derivation of the Hamiltonian  can be found in \cite{dorfler2013synchronization} (the authors call this an energy function).  The Hessian of the Hamiltonian can be expressed by:
\begin{align*}
  \nabla_{\tilde{x}\tilde{x}}\tilde{H}(\tilde{x},\tilde{u}) =
  \begin{bmatrix}
    \tilde{L}(\tilde{x}) & \\
    & A 
  \end{bmatrix}
\end{align*}
where $\tilde{L}(\tilde{x})\in\mathbb{R}^{(n-1)\times(n-1)}$ is defined as follows.
\begin{align*}
  \tilde{L}(\tilde{x})&:= 
  \begin{cases}
    \sum_{k\in\mathcal{N}(i)} B_{ik}\tilde{V}_i \tilde{V}_k\cos\theta_{ik} &\text{if}\quad i=j\\
    -B_{ij}\tilde{V}_i \tilde{V}_j\cos\theta_{ij} &\text{if}\quad i\neq j
  \end{cases}
\end{align*}
and $A$ is defined in \eqref{eqn:A}. An important observation is that $A$ is PD, and $\tilde{L}(\tilde{x})$ is PD if $\theta\in\Theta_G(\pi/2)$ because $\tilde{L}(\tilde{x})$ is a reduced Laplacian (see \cite{chung1997spectral,godsil2013algebraic,young2010robustness}). Thus, $\nabla_{\tilde{x}\tilde{x}}\tilde{H}(\tilde{x},\tilde{u})$ is PD if $\theta\in\Theta_G(\pi/2)$ and the Hamiltonian is locally strictly convex.

The Hamiltonian $\tilde{H}(\cdot,\cdot)$ can be obtained from that of the original port-Hamiltonian model \eqref{eqn:ss} by fixing the voltage magnitudes (i.e., $V_i\equiv\tilde{V}_i$ for $i\in\mathcal{V}$) and dropping the constants. Furthermore, the Hessian of $\tilde{H}(\cdot,\tilde{u})$ can be obtained by projecting the original Hessian on the space of $(\theta,\tomega)$. Consequently, we can see that the stability analysis based on a decoupling assumption is a special case of \eqref{eqn:ss}. Even though $\tilde{H}(\cdot,\tilde{u})$ is strictly convex on a reasonably large region, this can only be achieved if voltage magnitudes are assumed fixed.

\subsection{DC Approximation}
Under a DC approximation, where small angle differences are assumed, we replace $\sin\theta_{ij}$ with $\theta_{ij}$ and obtain
\begin{align*}
  \dot{\theta_i} &= \tomega_i - \tomega_n
  ,\quad i\in\mathcal{V}\setminus \{n\}\\
  \tau_{i} \dot{\tomega}_i&= -\tomega_i + \tomega^d - k_{P_i} (P_i - P_i^d)
  ,\quad i\in\mathcal{V}\\
  P_{i} &= \sum_{j\in\mathcal{N}(i)} B_{ij}\tilde{V}_i \tilde{V}_j \theta_{ij}, \quad i\in\mathcal{V}.
\end{align*}
The Hamiltonian is:
\begin{align*}
  \overline{H}(\tilde{x},\tilde{u}):=&\sum^n_{i=1} \left( \frac{\tau_{i}{\tomega}_i^2}{2k_{P_i}} 
  - P^u_i\theta_i\right)
  - \sum_{\{i,j\}\in\mathcal{E}} B_{ij} \tilde{V}_i \tilde{V}_j\frac{1}{2}\theta_{ij}^2.
\end{align*}
and we have that $ \dot{\tilde{x}} = (\tilde{R}-\tilde{J})\nabla_{\tilde{x}} \overline{H}(\tilde{x},\tilde{u})$. We thus have that condition \eqref{eqn:resilience} is satisfied and the Hessian is:
\begin{align*}
  \nabla_{\tilde{x}\tilde{x}}\overline{H}(\tilde{x},\tilde{u}) =
  \begin{bmatrix}
    \overline{L} & \\
    & A 
  \end{bmatrix}
\end{align*}
where $\overline{L}\in\mathbb{R}^{(n-1)\times(n-1)}$ is given by:
\begin{align*}
  (\overline{L})_{i,j} &= 
  \begin{cases}
    \sum_{k\in\mathcal{N}(i)} B_{ik}\tilde{V}_i \tilde{V}_k &\text{if}\quad i=j\\
    -B_{ij}\tilde{V}_i \tilde{V}_j &\text{if}\quad i\neq j
  \end{cases}
\end{align*}
and $A$ is defined in \eqref{eqn:A}. Matrix $\overline{L}$ is PD since it is reduced Laplacian and $A$ is also PD. As a result,$\nabla_{xx}\overline{H}(\tilde{x},\tilde{u})$ is PD everywhere. We thus have that the DC approximation yields a Hamiltonian that is {\em globally} strictly convex (the decoupling approximation only has a locally strictly convex Hamiltonian). The DC approximation, however, only holds when the voltage magnitudes are assumed fixed and angle differences are small.

\section{Computing Economic-Optimal and \\Stable Equilibria}\label{sec:economic}
In this section, we present optimization formulations to compute economic-optimal and stable equilibria for general \eqref{eqn:ss-gen} and port-Hamiltonian \eqref{eqn:ss} microgrid models. We first review classical OPF formulations and make connections with dynamic microgrid models. We then present a bilevel optimization formulation that is specialized to the port-Hamiltonian model \eqref{eqn:ss} and a probing formulation that can be applied to the general dynamic model \eqref{eqn:ss-gen}.

\subsection{Optimal Power Flow Formulation}
A typical form of OPF problems can be expressed by, but not limited to, the following form.
\begin{subequations}\label{eqn:classic}
  \begin{align}
    &\hspace{-.4in}\min_{\theta,V} \;  \sum_{i\in\mathcal{V}_G} c_i^1P_i + c_i^2(P_i)^2\\
    \st\; 
    &V^{l}_i\leq V_i \leq V^{U}_i,\; i\in\mathcal{V},\;
    |\theta_{ij}|\leq \theta^{U}_{ij},\; \{i,j\}\in\mathcal{E}\label{eqn:classic-x-1}\\
    &P^{L}_i\leq P_i \leq P^{U}_i,\;Q^{L}_i\leq Q_i \leq Q^{U}_i,\; i\in\mathcal{V}_G\label{eqn:classic-x-2}\\
    &P_i^d = P_i,\;Q^d_i = Q_i,\; i\in\mathcal{V}_L
  \end{align}
The readers are referred to \cite{molzahn2019survey,coffrin2018convex,coffrin2016qc} for more discussions on OPF problem formulations. Observe that introducing new variables $P^d_i$, $Q^d_i$, $\omega_i$, $V^d_i$ for $i\in\mathcal{V}_G$ and $\omega^d$ and enforcing the following constraints does not change the solution of \eqref{eqn:classic}.
\begin{align}
  0&= -\omega_i + \omega^d - k_{P_i} (P_i - P_i^d)\label{eqn:freq}
  ,&& i\in\mathcal{V}_G\\
  0&= -V_i + V_i^d - k_{Q_i} (Q_i - Q_i^d)
  ,&& i\in\mathcal{V}_G
\end{align}
\end{subequations}
Furthermore, defining $P^u_i$ and $Q^u_i$ as \eqref{eqn:u} and enforcing constraints $P^u_i=P^d_i$ and $Q^u_i=Q^d_i$ for $i\in\mathcal{V}_L$ does not change the solution. The following compact representation of OPFs are obtained:
\begin{align}\label{eqn:opf}
  \min_{x\in\mathbb{X},u\in\mathbb{U}}\; &c(x)\quad \st\;  f(x,u) = 0,
\end{align}
where 
\begin{align*}
  &c(x) := \sum_{i\in\mathcal{V}_G} c_i^1P_i^d + c_i^2(P_i^d)^2,\;\mathbb{X}:=\{x\mid \eqref{eqn:classic-x-1}-\eqref{eqn:classic-x-2}\text{ hold}\}\\
  &\mathbb{U}:=\{u\mid P^u_i = P^d_i,\; Q^u_i = Q^d_i,\quad\forall i\in\mathcal{V}_L\}
\end{align*}
and $f(\cdot,\cdot)$ is defined in \eqref{eqn:gen} (equivalently, one can write $g(x,u)=0$ and $h(x,u)=0$).
Thus, OPF problem can be regarded as finding an optimal steady-state. As discussed in the previous sections, the dynamic model \eqref{eqn:gen} does not have stability guarantee at the steady-state, and thus, the steady-state delivered by OPF is not guaranteed to be reachable by droop control. 
        
\subsection{Bilevel OPF for Port-Hamiltonian Models}
{The limitations of OPF motivate directly enforcing the stability in the optimization problem. Note that first solving the optimal power flow and then checking the stability (e.g., by using DAE simulation or by using Proposition \ref{prop:stab} for the port-Hamiltonian model) is not an option because if the stability condition is not satisfied with OPF solution, one needs to resort to a conservative, non-optimization-based operational decisions, which eventually translates to the loss of economic gain. Thus, optimization and stability enforcement should be performed simultaneously. Specifically, the goal is to minimize the given economic objective while satisfying the stability condition. To achieve such a goal, under the port-Hamiltonian assumption (Assumption \ref{ass:pHam}), one can consider formulating the problem as the following bilevel optimization problem:}
\begin{align}\label{eqn:blp}
  \min_{x\in\mathbb{X},u\in\mathbb{U}}\; &c(x) \quad \st\; x\in \amin_{\xi} H(\xi,u).
\end{align}
    {Here, $\amin(\cdot)$ is used to represent the set of local strict minima (with second-order sufficient conditions), and a dummy variable $\xi\in\mathbb{R}^{2n-1}\times\mathbb{R}^{n}_{>0}$ is used for the decision variables of the inner minimization problem. The constraint of \eqref{eqn:blp} can also be expressed as $\nabla_x H(x,u)=0$ and $\nabla_{xx}H(x,u)>0$. This formulation aims to achieve the dual goal of optimizing economic performance (by minimizing $c(\cdot)$ in the outer problem) and enforcing stability (by minimizing Hamiltonian). By Proposition \ref{prop:stab}, The solution of \eqref{eqn:blp} is {\em guaranteed} to be stable if Assumption \ref{ass:pHam} is satisfied.  The input $u$ is allowed to be varied in the outer problem so that one can seek for the best economic decision. 

To solve a bilevel problem, the inner problem is typically replaced by its stationarity conditions $\nabla_xH(x,u)=0$. We can easily show that such a problem is equivalent to OPF \eqref{eqn:classic} since $\nabla_x H(x,u)=0$ is equivalent to $f(x,u)=0$. In this sense, the  OPF problem \eqref{eqn:opf} can be seen as an approximation of the bilevel problem. However, there are currently no scalable methods for enforcing the second-order condition $\nabla_{xx}H(x,u)>0$ of the inner problem as constraints. Consequently, while the bilevel formulation achieves the desired stability properties of the solution, it is not computationally practical. Furthermore, the stability guarantee is fragile to different settings and network properties; if Assumption \ref{ass:pHam} does not hold, the solution is not guaranteed to be stable. {In other words, \eqref{eqn:blp} is only valid for port-Hamiltonian models.}

\subsection{Probing Formulation for General Models}\label{sec:probing}
{The computational challenge and limitations in flexibility of bilevel formulation motivate enforcing the stability by directly incorporating the dynamics into the optimization problem.} The following proposition implies that stability can be enforced by constraining probing trajectories around the equilibrium.

\begin{prop}\label{prop:lim}
  Consider system  $ \dot{x} =f(x,u^*)$ on $ \mathbb{R}^{n_x}$ and assume that there exists an equilibrium $x^*\in \mathbb{R}^{n_x}$ and $\delta>0$ such that, for any $x^0\in B_\delta(x^*)$, we have $\lim_{t\rightarrow\infty} x(t;x^0) = x^*$, where $B_\delta (x^*):=\{x\in\mathbb{R}^{n_x}\mid \Vert x-x^*\Vert<\delta \}$; then, $x^*$ is a locally asymptotically stable equilibrium.
\end{prop}
\begin{IEEEproof}
  The proof is trivial by definition of local asymptotic stability.
\end{IEEEproof}
Proposition \ref{prop:lim} suggests that we can find an economic-optimal and reachable equilibrium by solving the problem:
\begin{subequations}\label{eqn:sufprob}
  \begin{align}
    &\hspace{-0.4in}\inf_{x^*\in\mathbb{X},u^*\in\mathbb{U}} \quad c(x^*)\\
    \st \quad&f(x^*,u^*)=0\label{eqn:stnry0}\\
    &\dot{x}(t;x^0) = f(x,u^*),\; t \in\mathbb{R}_{\geq 0},\;x^0\in B_\delta (x^*)\label{eqn:dynopt}\\
    & \lim_{t\rightarrow\infty} x(t;x^0) = x^* ,\quad x^0\in B_\delta (x^*)    \label{eqn:convergence}
  \end{align}
\end{subequations}
Here, $x^*$ and $u^*$ are the state and input vector at the steady-state, $x(\cdot;x^0)$ is dynamic the state trajectory, and $x^0$ is the initial condition, and $B_\delta(x^*)$ can be interpreted as the region of attraction.

Problem \eqref{eqn:sufprob} is a semi-infinite optimal control problem; the semi-infinite nature arises from the fact that the differential equations must hold for all $t\in\mathbb{R}_{\geq 0}$ and for all $x^0\in B_\delta(x^*)$. We create a finite-dimensional approximation of the problem by (i) selecting a finite set of sample perturbations $\delta_s\in B_\delta(0)$ for $s\in\mathcal{S}$, where $\mathcal{S}$ is the set of samples. This approach can be interpreted as a stochastic programming problem in which one {\em probes} the equilibrium point $x^*$ by introducing random perturbations around it \cite{nemirovski2006scenario}. Furthermore, (ii) the differential equation constraints \eqref{eqn:dynopt} can be handled by applying discretization schemes (implicit Euler or orthogonal collocation \cite{biegler2010nonlinear}) with finite collocation points $t\in\mathcal{T}$, where $\mathcal{T}$ is the set of collocation time points. Finally, (iii) we relax the convergence condition \eqref{eqn:convergence} by replacing $t\rightarrow\infty$ with a sufficiently large end time $T\in\mathbb{R}_{>0}$ and  allowing small errors $\varepsilon\in\mathbb{R}_{>0}$ at the end time. This gives rise to the problem:
\begin{subequations}\label{eqn:approx}
  \begin{align}
    &\hspace{-0.4in}\min_{x^*\in\mathbb{X}, u^*\in\mathbb{U},\{x_s(t)\}_{s\in\mathcal{S},t\in\mathcal{T}}} \quad  c(x^*)\\
    \st \quad     &f(x^*,u^*)=0\label{eqn:stnry}\\
    &\sum_{t'\in\mathcal{T}}a_{t,t'}x(t')= f(x_s(t),u^*),\; t \in\mathcal{T},\;s\in\mathcal{S}\label{eqn:dynopt2}\\
    & x_s(0) = x^* + \delta_s, \; \Vert x_s(T) - x^* \Vert \leq \varepsilon,\; s\in\mathcal{S} \label{eqn:dynopt2-2}
  \end{align}
\end{subequations}
where $a_{t,t'}$ are the collocation coefficients. Observe that the semi-infinite problem \eqref{eqn:sufprob} is now transformed into a finite-dimensional NLP. This problem is high-dimensional but is also highly sparse and structured (the only coupling between scenario blocks are the equilibrium state variables $x^*$ and inputs $u$). This enables the use of structured NLP solvers that exploit parallel computers \cite{kang2015nonlinear}. 

\begin{rmk}[Interpretation of Probing as Stability Constraints]
  Note that \eqref{eqn:dynopt2}-\eqref{eqn:dynopt2-2} can be regarded as {\em stability constraints} that enforce reachability of the equilibrium provided by OPF. Accordingly, the suggested probing OPF \eqref{eqn:approx} can be regarded as a modification of OPF with stability constraints. Thus, the proposed framework can be applied along with the existing NLP-based OPF frameworks, since it can be performed by simply adding a set of constraints.
\end{rmk}
}

\section{Case Study}\label{sec:cstudy}
We illustrate the developments using the IEEE 118-bus system available at PGlib v18.08 \cite{pglib}. The data include: graph structure $\mathcal{G}(\mathcal{V},\mathcal{E},w)$, generator/load information $\mathcal{V}_G$ and $\mathcal{V}_L$, the generation cost coefficients $c^1$ and $c^2$, the voltage upper and lower bound $V^L$ and $V^U$, the upper and lower bound of active and reactive power generations $P^L$, $P^U$, $Q^L$, and $Q^U$, the active and reactive power loads $P^d$ and $Q^d$, and admittance matrix $Y$. We use dynamic paramters: $k_{P_i}=10$, $k_{Q_i}=1$, and $\tau_{P_i}=\tau_{Q_i}=10^{-3}$. We compare equilibria obtained using classical OPF \eqref{eqn:opf} and using the probing OPF \eqref{eqn:approx}. We present the results with the general model \eqref{eqn:ss-gen} as well as the results with the port-Hamiltonian \eqref{eqn:ss} model. The port-Hamiltonian model is obtained by modifying the general model with $G_{ij}:=0$ for $i,j\in\mathcal{V}$, $k_{P_i}:=100$, $k_{Q_i}:=10$, for $i\in\mathcal{V}_L$, $\mathcal{V}_G:=\mathcal{V}$, and $\mathcal{V}_L:=\emptyset$. For the approximation and obtaining probing samples, we use parameters: $T=1.0$, $\varepsilon=10^{-3}$, and $\delta_s\sim N(0_n,0.1^2 I_n)$. We use an orthogonal collocation scheme for the discretization \eqref{eqn:dynopt2}. We implement the NLPs with the algebraic modeling language {\tt JuMP} \cite{dunning2017jump} and solve them with NLP solver {\tt Ipopt} \cite{Wachter2006implementation}.

Economic performance is assessed by comparing the economic objective values evaluated at solutions $x^*$. For the general model, stability is assessed by simulating the dynamical system and observing the convergence. For the simulation, the initial state is perturbed by $x^0\leftarrow x^*+\delta^{\text{sim}}_s$ with random perturbations $\delta^{\text{sim}}_s\sim N(0_n, 0.1^2 I_n)$ and with $u^*$. The system is simulated over sufficiently long time $T_{\text{sim}}=1.0$. Note that initial condition samples for simulation were drawn from the same distribution with the probing samples but the actual samples were different. The stability is determined by checking whether $x\rightarrow x^*$ or not. In the simulation study, the errors $\theta_i(t)-\theta^*_i$, $\tomega_i(t)-\tomega^*_i$, and $V_i(t)-V_i^*$ for each node $i\in\mathcal{V}$ are plotted to verify convergence to the equilibrium. For the port-Hamiltonian model, the stability is also assessed by checking the eigenvalues of the Hessian $\nabla_{xx}H(x^*,u^*)$ of Hamiltonian.

\begin{figure*}[t!]
  \centering
  \includegraphics[width=.49\textwidth]{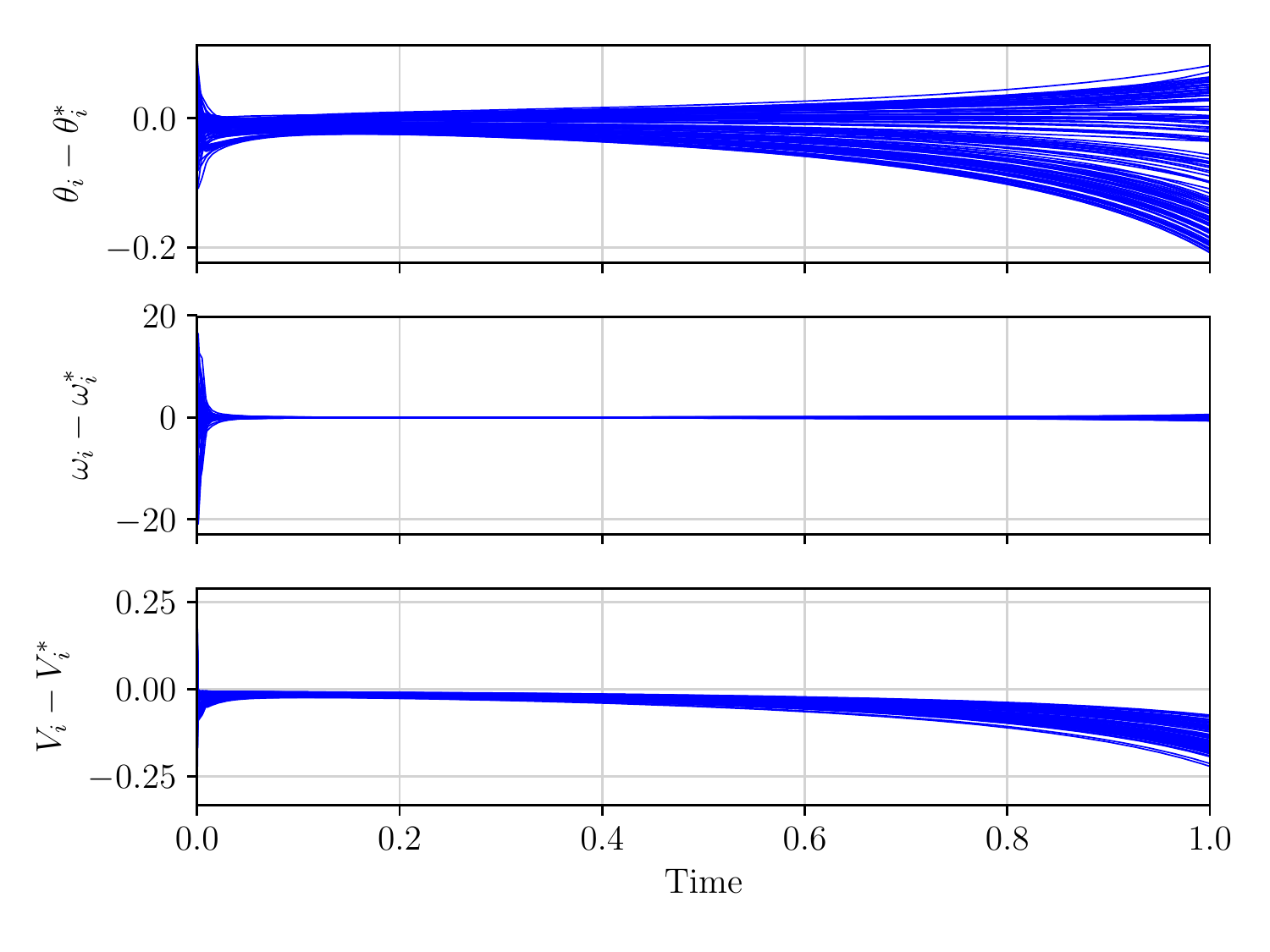}
  \includegraphics[width=.49\textwidth]{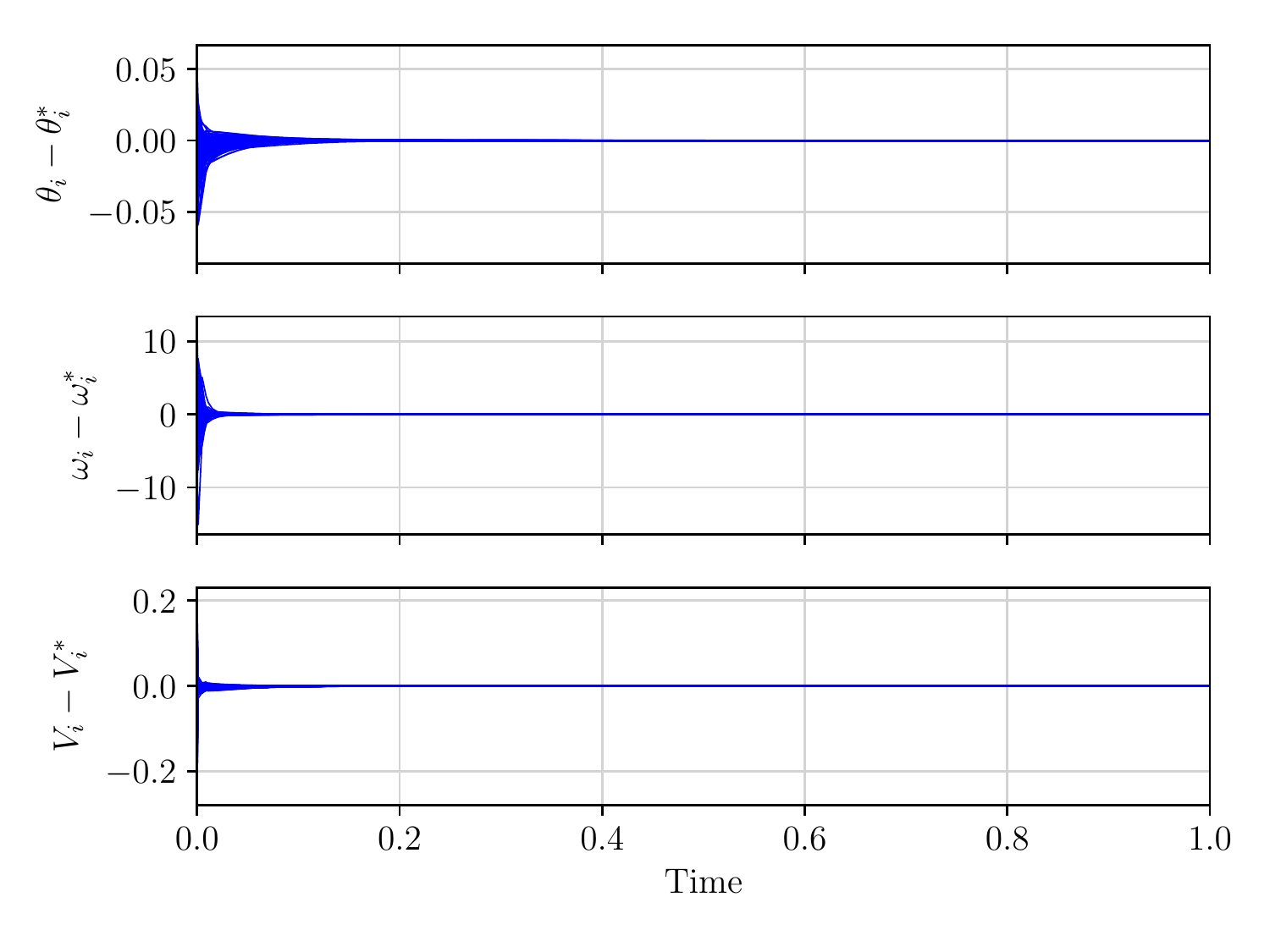}
  \caption{Simulation results for the {general model} with classical OPF (left) and the probing OPF with 4 scenarios (right).}\label{fig:sim-1}
  \centering
  \includegraphics[width=.49\textwidth]{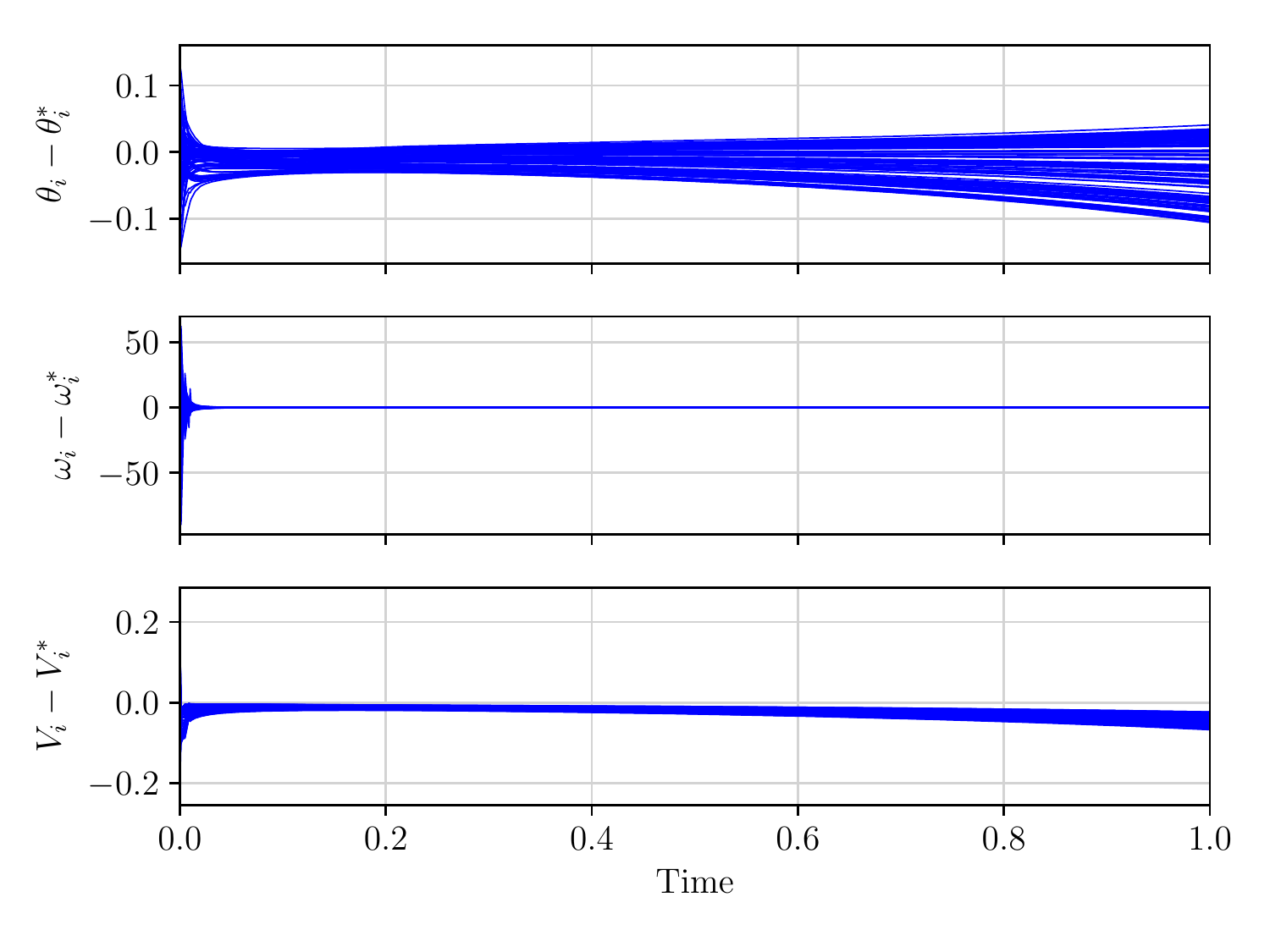}
  \includegraphics[width=.49\textwidth]{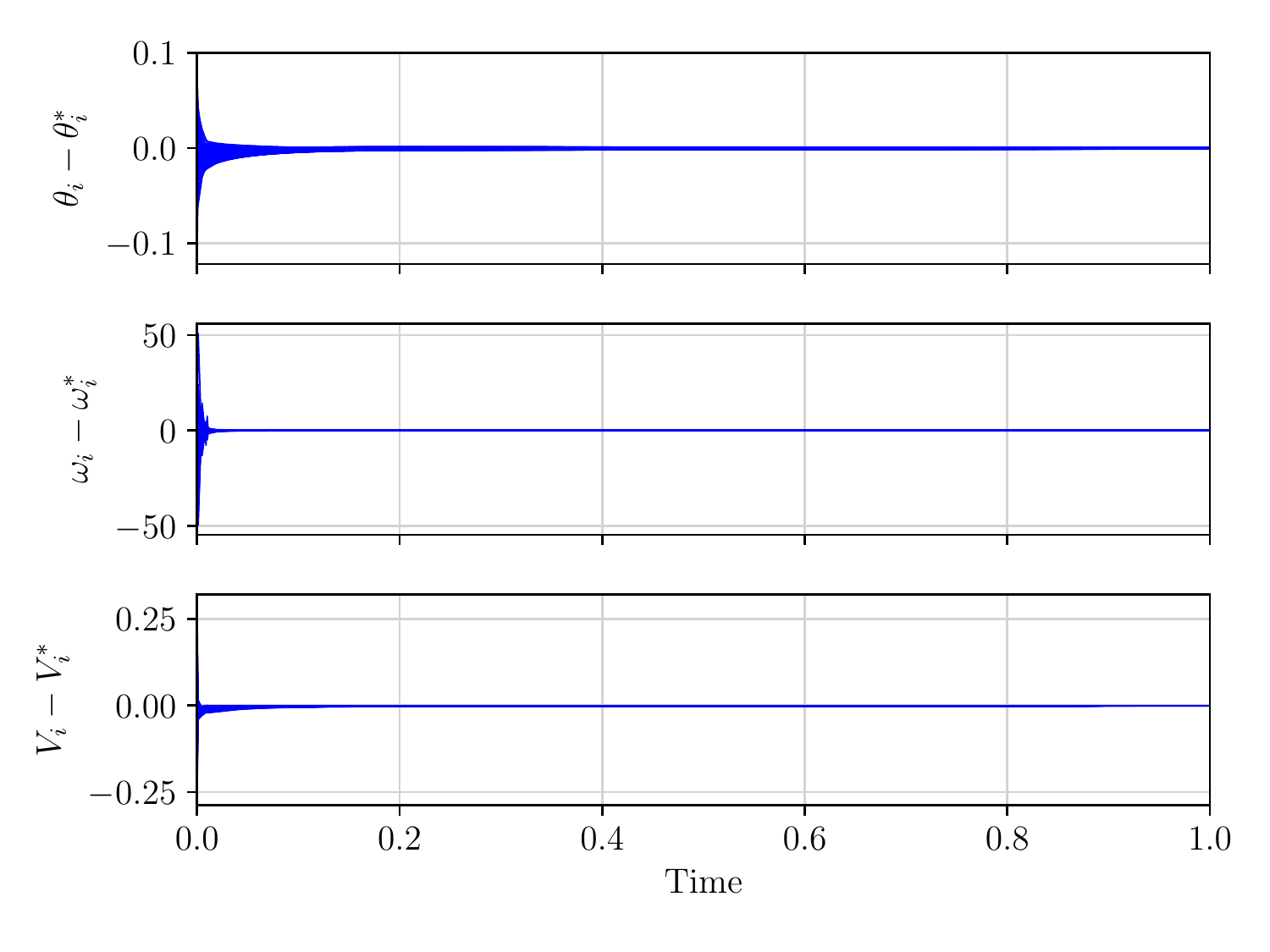}
  \caption{Simulation results for the port-Hamiltonian model with classical OPF (left) and the probing OPF with 4 scenarios (right).}\label{fig:sim-2}
\end{figure*}

The results for general and the port-Hamiltonian models are in Table \ref{table:res} and Fig \ref{fig:sim-1}-\ref{fig:sim-2}. The results reveal that for both general and port-Hamiltonian models, the equilibria identified with OPF achieve the best economic performances but they are not stable. Instability can be confirmed in the dynamic evolution shown in Fig \ref{fig:sim-1}-\ref{fig:sim-2}. In particular, for the port-Hamiltonian model, we can see that the Hessian has $1$ negative eigenvalues at the equilibrium (among a total of $353$ eigenvalues). This indicates that the equilibrium is a saddle point of Hamiltonian (thus, unstable).

\begin{table}[h!]
  \caption{Economics and stability results}
  \centering
  \label{table:res}
  \begin{tabular}{|c|c|c|c|c|}
    \firsthline
    &\multicolumn{2}{c|}{general model}&\multicolumn{2}{c|}{port-Hamiltonian model}\\
    \hline
    &$c(x^*)$  & stab.& $c(x^*)$ &$\lambda_{<0}(\nabla_{xx}H)$\\
    \hline
    OPF (0 sc.)& $1.136\times10^{5}$ & U/S&  $1.051\times10^{5}$ & 1\\
    \hline
    Pr. (1 sc.)& $1.138\times10^{5}$& U/S& $1.063\times10^{5}$& 1\\
    \hline
    Pr. (2 sc.)& $1.145\times10^{5}$& U/S& $1.075\times10^{5}$& 0\\
    \hline
    Pr. (4 sc.)& $1.175\times10^{5}$& S& $1.075\times10^{5}$& 0\\
    \lasthline
  \end{tabular}
  \vspace{0.1in}
  \caption{Computational statistics}
  \label{table:stat}
  \begin{tabular}{|c|c|c|c|c|c|}
    \firsthline
    \multicolumn{6}{|c|}{general model}\\
    \hline
    &\# var & \# eq. &\# ineq. &time (s) & \# iter.\\
    \hline
    OPF (0 sc.)&309 & 237& 358 & 0.266 &19\\
    \hline
    Pr. (1 sc.)&3495&3423&1066&5.341&59 \\
    \hline
    Pr. (2 sc.)&6681&6609&1774&8.603&39 \\
    \hline
    Pr. (4 sc.)&13053&12981&3190&32.489&48 \\
    \hline
    \multicolumn{6}{|c|}{port-Hamiltonian model}\\
    \hline
    &\# var & \# eq. &\# ineq. &time (s) & \# iter.\\
    \hline
    OPF (0 sc.)&309 & 237& 358 & 0.275 &26\\
    \hline
    Pr. (1 sc.)&2919&2837&938&2.603&34 \\
    \hline
    Pr. (2 sc.)&5529&5457&1518&6.873&34 \\
    \hline
    Pr. (4 sc.)&10749&10677&2678&22.796&44 \\
    \lasthline
  \end{tabular}
\end{table}

The results also show that the probing OPF is capable of finding a stable equilibrium for a small number of samples. For general and port-Hamiltonian models, when four samples are considered, the state trajectory converged to the steady-state $x^*$ (see Fig \ref{fig:sim-1}-\ref{fig:sim-2}). For the port-Hamiltonian model, this could also be confirmed by checking the number of negative eigenvalues (Table \ref{table:res}). This indicates that, as expected from Proposition \ref{prop:lim}, the dynamics encode significant stability information, so incorporating the probing trajectories can enforce stability. On the other hand, we note that the equilibria obtained with the probing formulation have inferior economic performance compared to OPF (the cost increases by 3.43\% for general model and 2.28\% for port-Hamiltonian, respectively). This indicates that there exists a strong trade-off between economics and stability. 

Computational statistics for the optimization problems are shown in Table \ref{table:stat}. We see that the probing solutions require more time than the OPF solution because they must capture the dynamics of the problem. Moreover, as expected, the solution time increases as the number of probing samples increase. We note, however, that the solution times are reasonable (in the range of less than one minute). For larger-scale problems, one can use parallel solvers to improve the solution times \cite{kang2015nonlinear}.

\section{Conclusions}\label{sec:conc}
We have presented a stochastic programming formulation for computing economic-optimal and stable equilibria for droop-controlled microgrids. Our work is motivated by the observation that standard OPF formulations might deliver economic-optimal set-points that are not reachable by droop control. We demonstrate that our approach is effective and can eliminate spurious (unreachable) equilibria obtained by standard OPF. Our work also reveals that, in general settings, an inherent trade-off between economic performance and stability might exist but that such a trade-off might disappear in certain restrictive settings (such as {constant voltage magnitude assumption}). {Our approach can also be used to identify economic-optimal stable equilibria for dynamical systems that arise in applications such as chemical reactors, robotics, and energy systems. }

\appendix

\subsection{Proof of Proposition \ref{prop:stab}}\label{sec:pf-stab}

We first claim that there exists a neighborhood $\mathbb{N}\subseteq \mathbb{D}$ of $x^*$ in which $x^*$ is the only stationary point (i.e. $\nabla_x H(x^*,u)=0$) of $H(\cdot,u)$ in $\mathbb{N}$. Suppose that such a neighborhood does not exist and let $\mathbb{N}_k := \{x\in\mathbb{D}\mid \Vert x-x^*\Vert < \frac{1}{k}\}$ for $k=1,2,\cdots$. Then there exists $x_k\in \mathbb{N}_k \setminus\{x^*\}$ such that $\nabla_x H(x_k,u)=0$. Let $y_k:=(x_k-x^*)/\Vert x_k-x^*\Vert$ for $k=1,2,\cdots$. Since each $y_k\in\mathbb{Y}:=\{y\in\mathbb{R}^{n_x}\mid \Vert y \Vert = 1\}$ and $\mathbb{Y}$ is compact, there exists a subsequence $\{y_{k_i}\}$ such that $\lim_{i\rightarrow\infty}y_{k_i}=y^*$ for some $y^*\in\mathbb{Y}$. Since the Hamiltonian is twice differentiable, we can use Taylor's theorem to establish that:
\begin{align}\label{eqn:Taylor}
  &\nabla_x H(x_{k_i},u) - \nabla_x H(x^*,u)\nonumber\\
  &= \nabla_{xx} H(x^*,u) (x_{k_i}-x^*) + R( x_{k_i} - x^*)
\end{align}
where the remainder term satisfies $\lim_{\Vert t\Vert \rightarrow 0}\frac{\Vert R(t)\Vert }{\Vert t\Vert}=0$.
Using $\nabla_x H(x_{k_i}) = 0$ for $i=1,2,\cdots$ and dividing \eqref{eqn:Taylor} by $\Vert x_{k_i}-x^*\Vert$, we have that:
\begin{align}
  \nabla_{xx} H(x^*,u) {y_{k_i}} = -\frac{R(x_{k_i}-x^*)}{\Vert x_{k_i}-x^*\Vert}.
\end{align}
By taking $i\rightarrow\infty$, we have that $\nabla_{xx}H(x^*,u)y^* = 0$.  Since $\nabla_{xx}H(x^*,u)$ is PD, $y^*=0$, but this contradicts $y^*\in\mathbb{Y}$. Therefore, there exists a neighborhood $\mathbb{N}$ of $x^*$ such that $x^*$ is the only stationary point in $\mathbb{N}$.

Let $\mathbb{N}$ be such a neighborhood. From \eqref{eqn:resilience2}, $\nexists\, x\in \mathbb{N}$ such that $\dot{H}(x,u)=0$. Therefore, $\dot{H}(x,u)<0$ for any $x\in \mathbb{N}$, we can construct a Lyapunov function $V:\mathbb{N}\rightarrow \mathbb{R}$ of the form $V(x):=H(x,u)-H(x^*,u)$. The function satisfies $\dot{V}(x)<0$ and $V(x)>0$ in $x\in \mathbb{N}\setminus\{x^*\}$ as well as $V(x^*)=0$.  By Lyapunov's stability Theorem \cite[Theorem 4.1]{khalil1996nonlinear}, $x^*$ is locally asymptotically stable.

\subsection{Proof of Proposition \ref{prop:convex}}\label{sec:pf-convex}
We observe that $\nabla_{xx}H(x,u)$ is PD if and only if $D(u)+T(x)$ and its Schur complement are PD \cite[Theorem 7.7.6]{horn1990matrix}. By noting that $A$ is PD, it suffices to show that $D(u)+T(x)$ and $S(x,u)$ (defined in \eqref{eqn:S}) are PD. $D(u)$ is PD from the  assumption that $Q^d_i+V^d_i/k_{Q_i}>0$. We consider $ T(x) = T_1(x) + T_2(x)$ where 
\begin{align*}
  (T_1(x))_{i,j} &=
  \begin{cases}
    -B_{ii}-\sum_{k\in\mathcal{N}(i)}B_{ik} \cos(\theta_i-\theta_k) &\quad\text{if}\quad i=j
    \\
    0& \quad\text{if}\quad i\neq j
  \end{cases}\\
  (T_2(x))_{i,j} &=
  \begin{cases}
    \sum_{k\in\mathcal{N}(i)} B_{ik} \cos(\theta_i-\theta_k) &\quad\text{if}\quad i=j 
    \\
    -B_{ij} \cos\theta_{ij}&\quad\text{if}\quad i\neq j
  \end{cases}
\end{align*}
By the definition of admittance matrix, $-B_{ii}=\sum_{k\in\mathcal{N}(i)}B_{ik}$ for $i\in\mathcal{V}$. This implies that $(T_1(x))_{i,i}=\sum_{k\in\mathcal{N}(i)} B_{ik}(1-\cos(\theta_i-\theta_k))\geq 0$ for $i\in\mathcal{V}$. Thus, $T_1(x)$ is PSD. Furthermore, $T_2(x)$ is PSD since it takes a weighted Laplacian form. Therefore, $D(u)+T(x)$ is PD and its smallest eigenvalue is greater than or equal to the smallest eigenvalue of $D(u)$. This yields that:
\begin{align}\label{eqn:D+T}
  \lambda_{\max}((D(u)+T(x))^{-1})\leq 1/\lambda_{\min}(D(u)).
\end{align}

Since $(D(u)+T(x))^{-1}$ is PD, $W(x)^\top(D(u)+T(x))^{-1} W(x)$ is PSD. We have that:
\begin{align*}
  &\lambda_{\max}(W(x)^\top(D(u)+T(x))^{-1} W(x))\nonumber\\
  &=\max_{\Vert v\Vert = 1}(W(x)v)^\top(D(u)+T(x))^{-1} (W(x)v)\nonumber\\
  &=\max_{ u=W(x) v,\Vert v\Vert = 1}u^\top(D(u)+T(x))^{-1} u\nonumber\\
  &\leq\max_{ \Vert u\Vert \leq \Vert \tilde{W}(x)\Vert }u^\top(D(u)+T(x))^{-1} u\nonumber\\
  &= \lambda_{\max}((D(u)+T(x))^{-1}) \Vert \tilde{W}(x) \Vert^2.
\end{align*}
Here, $\tilde{W}(x)\in\mathbb{R}^{n\times n}$ and $\tilde{W}(x)=[W(x)^\top \; 0_n]$, where $0_n\in\mathbb{R}^n$ is a zero vector.
By \cite[pg 314]{horn1990matrix}, $\Vert \tilde{W}(x) \Vert^2 \leq \Vert \tilde{W}(x) \Vert_1^2$ ($\Vert \cdot \Vert_1$ denotes $\ell 1$-norm), and $\Vert \tilde{W}(x) \Vert_1^2 = \Vert W(x) \Vert_1^2$. This yields:
\begin{align*}
  &\lambda_{\max}(W(x)^\top(D(u)+T(x))^{-1} W(x))
  \leq \frac{\Vert W(x) \Vert_1^2}{\lambda_{\min}(D(u))}.
\end{align*}
We also have that
\begin{align*}
  \Vert W(x)\Vert_1 &=\sum_{i\in\mathcal{V}\setminus \{n\}}\sum_{j\in\mathcal{V}} |(W(x))_{i,j}|\nonumber\\
  &\leq 3B_{\max}V_{\max}|\mathcal{E}| \epsilon
\end{align*}
where $B_{\max}:=\max\{|B_{ij}|\mid \{i,j\}\in\mathcal{E}\}$.  This yields:
\begin{align}\label{eqn:conv-1}
  &W(x)^\top (D(u)+T(x))^{-1}W(x)
  \nonumber\\
  &\leq \frac{9(B_{\max})^2(V_{\max})^2|\mathcal{E}|^2 \epsilon^2}{\lambda_{\min}(D(u))} I.
\end{align}
Let $\lambda_{\min}(L_{\pi/4})$ be the smallest eigenvalue of $L_{\pi/4}\in\mathbb{R}^{(n-1)\times(n-1)}$ where
\begin{align*}
  (L_{\pi/4})_{i,j} &= 
  \begin{cases}
    \sum_{k\in\mathcal{N}(i)} B_{ik}(V_{\min})^2 \cos(\pi/4) &\text{if}\quad i=j\\
    -B_{ij}(V_{\min})^2 \cos(\pi/4) &\text{if}\quad i\neq j.
  \end{cases}
\end{align*}
The smallest eigenvalue of $L(x)$ is greater than or equal to $\lambda_{\min}(L_{\pi/4})$ if $\epsilon\leq \pi/4$. This yields:
\begin{align}\label{eqn:conv-2}
  L(x)
  \geq \lambda_{\min}(L_{\pi/4})I.
\end{align}

By \eqref{eqn:conv-1} and \eqref{eqn:conv-2}, we have that
\begin{align*}
  S(x,u)
  \geq \left[\lambda_{\min}(L_{\pi/4})- \frac{9(B_{\max})^2(V_{\max})^2|\mathcal{E}|^2 \epsilon^2}{{\lambda_{\min}(D(u))}} \right]I.
\end{align*}
$\lambda_{\min}(L_{\pi/4})$, $B_{\max}$, $V_{\max}$, $\lambda_{\min}(D(u))$, $|\mathcal{E}|$ are independent of the choice of $x$. If $\epsilon$ is sufficiently small, then $S(x,u)$ is PD and so is $\nabla_{xx}H(x,u)$.

\ifCLASSOPTIONcaptionsoff
\newpage
\fi

\bibliographystyle{IEEEtran}
\bibliography{droop}
\begin{IEEEbiographynophoto}{Sungho Shin}
  received his B.S. in chemical engineering and mathematics from Seoul National University, South Korea in 2016. He is currently a Ph.D. candidate in the Department of Chemical and Biological Engineering at the University of Wisconsin-Madison. His research interests include control theory, optimization algorithms, and complex networks.
\end{IEEEbiographynophoto}
\vspace{-0.1in}

\begin{IEEEbiographynophoto}{Victor M. Zavala}
  is the Baldovin-DaPra Associate Professor in the Department of Chemical and Biological Engineering at the University of Wisconsin-Madison. He holds a B.Sc. degree from Universidad Iberoamericana and a Ph.D. degree from Carnegie Mellon University, both in chemical engineering. He is an associate editor for the Journal of Process Control and a technical editor of Mathematical Programming Computation. His research interests are in the areas of energy systems, high-performance computing, stochastic programming, and predictive control.
\end{IEEEbiographynophoto}

\vfill

\end{document}